\documentclass[11pt]{article}
\usepackage{sustyle}
\usepackage{tikz}
\usetikzlibrary{arrows.meta}
\usetikzlibrary{arrows}
\tikzstyle{main node}=[draw]

\usepackage[T1]{fontenc}
\usepackage{amsmath, graphicx, tikz, enumerate, amssymb, pgf}
\usepackage[colorlinks=true, allcolors=black,c itecolor=blue, urlcolor=black]{hyperref}
\usepackage{bibentry}
\usepackage{setspace}
\usepackage{pgfplots}
\usepackage{bchart}
\def\d{\mathrm{d}}
\def\ep{\epsilon}
\def\del{\delta}
\newcommand\dist{D}
\usepackage{enumitem}
\usepackage{booktabs}

\DeclareMathOperator{\TV}{TV}
\DeclareMathOperator{\id}{id}
\DeclareMathOperator{\kl}{kl}
\newcommand{\citep}[1]{\cite{#1}}
\newcommand{\citet}[1]{\cite{#1}}

\usepackage{mathabx}

\newtheorem{axiom}[theorem]{Axiom}

\title{A Statistical Viewpoint on Differential Privacy: Hypothesis Testing, Representation and Blackwell's Theorem}

\date{October 21, 2024}

\begin{document}

\author{Weijie J.\ Su\thanks{University of Pennsylvania. Email: \texttt{suw@wharton.upenn.edu}.}}
\maketitle



\begin{abstract}

Differential privacy is widely considered the formal privacy for privacy-preserving data analysis due to its robust and rigorous guarantees, with increasingly broad adoption in public services, academia, and industry. Despite originating in the cryptographic context, in this review paper we argue that, fundamentally, differential privacy can be considered a \textit{pure} statistical concept. By leveraging David Blackwell's informativeness theorem, our focus is to demonstrate based on prior work that all definitions of differential privacy can be formally motivated from a hypothesis testing perspective, thereby showing that hypothesis testing is not merely convenient but also the right language for reasoning about differential privacy. This insight leads to the definition of $f$-differential privacy, which extends other differential privacy definitions through a representation theorem. We review techniques that render $f$-differential privacy a unified framework for analyzing privacy bounds in data analysis and machine learning. Applications of this differential privacy definition to private deep learning, private convex optimization, shuffled mechanisms, and U.S.\ Census data are discussed to highlight the benefits of analyzing privacy bounds under this framework compared to existing alternatives.

\end{abstract}

\section{Introduction}

In the era of big data, the increasing collection and analysis of personal information have raised significant concerns about privacy. It is therefore imperative to develop general approaches capable of conducting large-scale data analysis while still preserving the privacy of the individuals whose information is included in the data. Differential privacy, introduced by \cite{DMNS06}, has emerged as a rigorous mathematical framework for quantifying and reasoning about privacy risks associated with data analysis.

Owing to many of its appealing features \citep{dwork2014algorithmic}, this privacy notion has gained significant attention in both academia and industry \citep{cummings2024advancing}, and has been adopted by major technology companies including Google \citep{rappor,aktay2020google,xu2023federated}, Apple \citep{apple}, Microsoft \citep{microsoft}, and LinkedIn \citep{rogers2020linkedin}. In 2020, the U.S.\ decennial census adopted differential privacy when it released demographic data \citep{abowd2018us,abowd20222020}.\footnote{For an overview of real-world deployments of differential privacy and related software development, see \url{https://desfontain.es/blog/real-world-differential-privacy.html}.}

Despite being proposed by the cryptographic community, differential privacy was recognized early on to have a strong connection to statistics \citep{wasserman_zhou}. Roughly speaking, differential privacy's mental model of distinguishing whether an individual is present or absent can be formulated as a hypothesis testing problem. All existing differential privacy definitions leverage this concept in their formulations, some more directly than others \citep{dwork2014algorithmic,approxdp,concentrated,concentrated2,tcdp,renyi}. However, it has long been unclear whether hypothesis testing is the right language for reasoning about differential privacy in a rigorous sense and, if so, how we can obtain a unified treatment of all existing differential privacy definitions to enhance the framework by incorporating statistical insights and techniques.

In \citet{KOV,dong2022gaussian}, the authors formally addressed this question in the affirmative. The core of their strategies is to leverage a result due to David Blackwell to identify a measure that encodes the complete information for hypothesis testing in the context of differential privacy \citep{blackwell1950comparison}. Roughly speaking, this information-complete measure is given by the envelope of all uniformly most powerful tests using the likelihood ratio, as characterized by the Neyman--Pearson lemma. With this measure in place, \citet{dong2022gaussian} introduced $f$-differential privacy and showed, using a representation theorem, that any differential privacy definition could be considered a subfamily of this privacy definition.

The definition of $f$-differential privacy has many appealing features, owing to its tightness inherent to its definition \citep{dong2021central,bu2020deep}. For example, it is tight and information-lossless under composition and subsampling, which are fundamental operations in private machine learning. Since it was proposed, there has been an increasing line of research focusing on leveraging $f$-differential privacy in applications where optimal privacy-utility trade-off is desired, including private deep learning, private convex optimization, census privacy, and privacy auditing~\citep{bu2020deep,bok2024shifted,su2024revealing,nasr2023tight}.

In this review, we aim to provide a minimal elaboration on the $f$-differential privacy framework as well as a selective survey of recent results using $f$-differential privacy for enhanced privacy analysis of important privacy-preserving machine learning methods, with an emphasis on highlighting the statistical insights behind the development. We remark that the literature is quickly growing and this review is by no means exhaustive. In Section~\ref{sec:basics-dp-without}, we start with the fundamentals of differential privacy to elucidate why $f$-differential privacy is the natural definition from the statistical viewpoint. In Section~\ref{sec:privacy-accounting}, we introduce the basics of the framework of $f$-differential privacy and useful techniques accompanying this framework. In Section~\ref{sec:applications}, we showcase the effectiveness of $f$-differential privacy in obtaining improved or even optimal privacy analysis for several important applications. In Section~\ref{sec:future-directions}, we identify challenges and future research directions for differential privacy from a statistical viewpoint.


\section{Differential Privacy from a Statistical Viewpoint}\label{sec:basics-dp-without}

\subsection{Background}
\label{sec:background}

Consider a setting where a trusted curator holds a dataset with each row containing an individual's data \citep{dwork2006differential}. A (private) algorithm or analysis, depending on the context, takes as input the dataset and produces an output. The goal is to protect the privacy of every individual while allowing meaningful statistical analysis of the dataset in the sense that the output is approximately close to that of the non-private counterpart of the algorithm.

Various attacks based on the output of the algorithm can be launched against a dataset to compromise individual privacy \citep{dwork2017exposed}. Reconstruction attacks aim to recreate the entire dataset or a significant portion of it, potentially exposing sensitive information about individuals~\citep{dinur2003revealing}. Tracing attacks~\citep{homer2008resolving}, or membership inference attacks in the setting of access to black-box machine learning models~\citep{shokri2017membership}, attempt to determine the presence or absence of a specific individual's record in the dataset.

These attacks highlight the need for a rigorous and principled approach to privacy protection, as defending against each attack individually would be endless and ad hoc. For this purpose, it is beneficial to recognize that a common feature of the attacks above is that the success of an attack essentially allows an adversary to determine that some reality is true instead of the other. For instance, a successful tracing attack could reveal whether a specific individual is present in the dataset or not.

In light of this,  \cite{DMNS06} envisioned the problem of determining which of two possible realities is real, where the two possibilities correspond to \textit{neighboring} datasets $S$ and $S'$ in the sense of differing in only one individual's record. That is, $S$ and $S'$ are identical except for one individual such that $S$ has, for example, Alice's data while $S'$ has Anne's data. If distinguishing between these two realities is fundamentally difficult, regardless of which individual's record differs, then any attack that violates an individual's privacy would necessarily fail.

The notion of differential privacy~\citep{DMNS06} well captures the essence of imposing the difficulty in distinguishing the two possibilities. The name ``differential'' underscores its emphasis on the fact that the two competing realities or, explicitly, datasets, differ only in one individual's data. By satisfying differential privacy, an algorithm can resist all possible attacks that ultimately lead to privacy violation of some individual. This allows researchers to focus on designing algorithms that provably meet the privacy guarantee.

Roughly speaking, differential privacy limits an attacker's ability to discern the presence or absence of a single record while still allowing for the discovery of useful statistical patterns that do not rely on any one individual's data. This property of differential privacy enables the extraction of valuable insights from datasets without compromising the privacy of the individuals contained within them. For example, releasing aggregate statistics that reveal correlations between variables, say, smoking and lung cancer, does not violate differential privacy, even if the statistics can be used to predict sensitive information about an individual, as long as the discovery of these findings is not dependent on any single individual's data~\citep{dwork2017guilt,ullman2021statistical}.

\subsection{Hypothesis Testing Interpretation}
\label{sec:metrics}

Having introduced the background that motivated differential privacy, in this section we formally define this privacy notion. The ``differential'' nature of differential privacy can be precisely formulated as a hypothesis testing problem in which an attacker makes a choice between
\begin{equation}\label{eq:ht}
\begin{aligned}
H_0&:  \text{algorithm } M \text{ is applied to dataset } S \\
H_1&:  \text{algorithm } M \text{ is applied to dataset } S'.
\end{aligned}
\end{equation}
This connection was first recognized by \cite{wasserman_zhou} and further refined by \cite{KOV}. To interpret this hypothesis testing problem, we envision a powerful attacker who knows both the output of the algorithm and its probability distributions under both the null $H_0$ and alternative $H_1$. This probability distribution fully characterizes the information provided by the algorithm. By assuming such a strong attacker, we ensure that the resulting privacy definition, which can render the attacker powerless, is inherently strong.

Having formulated differential privacy from the hypothesis testing perspective, the next step is to specify the basis upon which the attacker determines which hypothesis to reject or accept. This necessitates a metric between the two probability distributions corresponding to the null and the alternative to quantify their difference. This motivates us to propose the following axiom as a
self-evident requirement for any possible differential privacy definition.

\begin{axiom}\label{axiom:ht}
\textit{Any differential privacy definition should be defined via hypothesis testing through some metric between the probability distributions of the output when the algorithm is applied to two neighboring datasets.}
\end{axiom}

We elaborate on this axiom as follows. Consider a new differential privacy definition using a metric we call $\dist$.\footnote{Here, we use the term ``metric'' in a broader sense, as $\dist$ may not necessarily be univariate, symmetric, or satisfy the triangle inequality.} The difference between the probability distributions in Equation~\ref{eq:ht} can be quantified by $\dist(M(S), M(S'))$, where we identify $M(S)$ with the probability distribution of the output of $M$ when applied to the dataset $S$. The metric takes values in the parameter space ${\xi}$ endowed with a partial order. For two privacy parameters $\xi_1$ and $\xi_2$, if $\xi_1 \preceq \xi_2$ or equivalently $\xi_2 \succeq \xi_1$, we say that $\xi_1$ delivers a stronger privacy guarantee than $\xi_2$ does in the sense of this differential privacy definition.

Formally, we say that $M$ has $\xi$-differential privacy with some parameter $\xi$ if
\[
\dist(M(S), M(S')) \preceq \xi
\]
for any neighboring $S$ and $S'$. As we shall soon see, the literature has used the choices of the supremum of log-likelihood ratio, R\'enyi divergence, and others to define differential privacy \citep{DMNS06,approxdp,concentrated,concentrated2,tcdp,renyi,dong2022gaussian}.

Axiom~\ref{axiom:ht} alone does not fully capture the essence of the hypothesis testing interpretation of differential privacy. Another crucial aspect to consider is the \textit{post-processing} property. This property stems from the intuitive notion that the privacy level of a statistical analysis should not be compromised by simply processing the output of the algorithm $M$. In essence, any meaningful definition of differential privacy must inherently satisfy this post-processing property.

We consider post-processing a basic requirement for any meaningful differential privacy definition. It is so fundamental that we consider it an axiom.
\begin{axiom}\label{axim}
\textit{Any differential privacy definition must satisfy the property of post-processing.}
\end{axiom}

Formally, denote by $\mathrm{Proc}$ any post-processing operation on the output of the algorithm, such as dimension reduction or summarization. Axiom~\ref{axim} states that 
\[
\dist(M(S), M(S')) \text{ and } \dist(\mathrm{Proc}(M(S)), \mathrm{Proc}(M(S')))
\] 
must be comparable and satisfy
\[
\dist(\mathrm{Proc}(M(S)), \mathrm{Proc}(M(S'))) \preceq \dist(M(S), M(S')).
\]
Notably, this requirement precludes Wasserstein distances because this inequality does not hold when $\mathrm{Proc}$ is an affine transformation with a scale factor greater than 1.

\subsection{A Representation Theorem}
\label{sec:repr-theor}

Taken together, we show that Axioms~\ref{axiom:ht} and \ref{axim} require the metric to have a certain structure to make sense for defining differential privacy. To state this result, we define the trade-off function for a pair of probability distributions defined on the same space. Before doing so, we need to formally define the trade-off function between $P$ and $Q$, which can be defined through the likelihood ratio using the Radon--Nikodym derivative. For any constant $r$, consider the pair
\begin{equation}\label{eq:lration}
\alpha := \P_{P}\left( \frac{\d P}{\d Q} \le r \right), \quad \beta := \P_{Q}\left( \frac{\d P}{\d Q} > r \right).
\end{equation}
By varying the value of $r$, this parametric equation gives rise to the trade-off function, denoted as $T(P, Q): [0, 1] \longmapsto [0, 1]$, which is defined as\footnote{Due to possible discontinuity of the likelihood function, we add linear interpolation if some points disconnect. Alternatively, the final trade-off function is given by the largest convex lower envelope of $T(P, Q)$ defined on all ``continuous'' $\alpha$ in Equation \ref{eq:lration}.}
\[
T(P, Q)(\alpha) = \beta.
\]

\begin{theorem}[Representation Theorem]\label{thm:full}
Under Axioms~\ref{axiom:ht} and \ref{axim}, any differential privacy definition must have its metric $\dist$ depend on the probability distributions through the trade-off function. That is, there must exist a link function $d$ defined on trade-off functions such that
\begin{equation}\label{eq:represent}
\dist(P, Q) = d(T(P, Q)).
\end{equation}
\end{theorem}

\begin{remark}
A proof of Theorem~\ref{thm:full} is given in Section~\ref{sec:blackw-persp}. For the construction of link functions for certain types of metrics $D$, see Proposition B.4 in \citet{dong2022gaussian}. This theorem in the special case of $(\epsilon, \delta)$-differential privacy was first obtained by \citet{KOV}.
\end{remark}

This theorem demonstrates that any differential privacy definition admits a representation through the trade-off function. Roughly speaking, any differential privacy definition exploits no more information than conveyed by the trade-off function in an information-theoretic sense. Recognizing this fact, \citet{dong2022gaussian} considered directly using the trade-off function for definition by setting $d$ to the identity, which leads to $f$-differential privacy. Letting $f$ denote a general trade-off function, an algorithm $M$ is $f$-differentially private if\footnote{Here, we write $g \ge_\preceq f$ or, equivalently, $f \le_\succeq g$ if $g(\alpha) \ge f(\alpha)$ for all $0 \le \alpha \le 1$. The design of the notation $\ge_\preceq$ is to reflect that $g$ is numerically greater than $f$; meanwhile, the two distributions associated with $g$ are closer to each other than those associated with $f$.}
\[
T(M(S), M(S')) \ge_{\preceq} f
\]
for any neighboring datasets $S$ and $S'$. Theorem~\ref{thm:full} shows that any other differential privacy definition is necessarily an information-lossy reduction of $f$-differential privacy.


Now we list important differential privacy variants and show how the representation theorem relates to each of them. This is not to say that any link function $d$ suffices for a new definition. Theorem~\ref{thm:full} provides a necessary instead of a sufficient condition for a reasonable differential privacy definition. For example, the total variation distance, which corresponds to
\[
d_{\TV}(f) = \frac12 \int_0^1 |1 + f'(x)| \d x,
\]
is not inappropriate for a privacy definition~\citep{DMNS06,barber2014privacy}. See \citet{desfontaines2019sok} for a comprehensive overview of differential privacy definitions.

\begin{enumerate}
\item[] {\bf Pure differential privacy} 

The max-divergence
\begin{equation}\nonumber
D_{\infty}(P, Q) := \max \log \frac{\d P}{\d Q}
\end{equation}
is used to define pure differential privacy \citep{DMNS06}. The corresponding link function is $d_{\infty}(f) = -\log(-f'(1))$, where the derivative is understood as the left derivative. It is possible that $f'(1) = 0-$. Thus, $\ep = \infty$ is included for completeness.

When both $P$ and $Q$ have probability density functions or probability mass functions, denoted $p(x)$ and $q(x)$, respectively, then $D_{\infty}(P, Q)$ equals the supremum of $\log \frac{p(x)}{q(x)}$. An algorithm $M$ is $\ep$-differentially private if $D_{\infty}(M(S), M(S'))$ is upper bounded by $\ep \ge 0$ for any neighboring datasets $S$ and $S'$.

\item[] {\bf Approximate differential privacy} 

For approximate differential privacy \citep{approxdp}, we fix $0 \le \del < 1$ and consider $\ep$ the only privacy parameter:
\[
D_{\infty}^{\delta}(P, Q) = \max_{\Omega: P(\Omega) \ge \delta} \log \frac{P(\Omega) - \delta}{Q(\Omega)} = \inf_{\ep} \{\ep: H_{\e^{\ep}}(P, Q) \le \delta\},
\]
where $H_{\e^{\ep}}(P, Q)$ is the Hockey-stick divergence defined as 
\[
H_{\e^{\ep}}(P, Q) = \E_{Q} \left( \frac{\d P}{\d Q} - \e^{\ep} \right)_+,
\] 
where $x_+$ denotes $\max\{x, 0\}$. An algorithm $M$ satisfies $(\epsilon, \delta)$-differential privacy if 
\[
D_{\infty}^{\delta}(M(S), M(S'))
\]
is upper bounded by $\epsilon$ for any neighboring datasets $S$ and $S'$. Smaller values of $\epsilon$ and $\delta$ provide stronger privacy protection. More precisely, $\delta$ should be set to a very small value, while a relatively large value of $\epsilon$ can still provide meaningful privacy guarantees \citep{approxdp,near2022differential}.

Regarding Theorem~\ref{thm:full}, the link function $d_{\infty}^{\delta}$ at $f$ can be defined as the unique root in $\ep \ge 0$ of the following equation \citep{dong2022gaussian}:\footnote{If the root does not exist, which could be the case if $\delta$ is small, $\epsilon$ is set to $\infty$ as a convention.}
\begin{equation}\label{eq:dual}
f^\ast(-\e^{-\ep}) + \e^{-\ep}(1 - \delta) = 0,
\end{equation}
where $f^\ast(x) = \sup_y (xy - f(y))$ is the convex conjugate of $f = T(P, Q)$. This representation is implied by a connection between $(\epsilon, \delta)$-differential privacy and hypothesis testing \citep{wasserman_zhou,KOV}. By varying $\del$, all $(\ep, \del)$-differential privacy guarantees taken together give the privacy profiles of the algorithm \citep{balle2018privacy}, which are the dual representation of the trade-off function $f$.

\item[] {\bf Divergence-based differential privacy} 

The concentrated differential privacy family of definitions \citep{concentrated,concentrated2,tcdp} and R\'enyi differential privacy \citep{renyi} are defined based on the R\'enyi divergence
\[
\dist_{\gamma}(P, Q) := \frac1{\gamma-1} \log \E_{Q} \left(\frac{\d P }{\d Q} \right)^\gamma
\]
for $\gamma > 1$. For example, an algorithm $M$ is said to be $(\gamma, \ep)$-R\'enyi differentially private if
\[
\dist_{\gamma}(M(S), M(S')) \le \ep.
\]

Taking $f = T(P, Q)$, Equation~\ref{eq:represent} in Theorem~\ref{thm:full} for R\'enyi differential privacy takes the form
\[
D_{\gamma}(P, Q) = d_{\gamma}(f) = \frac1{\gamma-1} \log \int_0^1 |f'(x)|^{1-\gamma} \d x.
\]

\end{enumerate}

\subsection{The Blackwell Perspective}
\label{sec:blackw-persp}

The representation theorem can be proved using a celebrated result due to David Blackwell.

\begin{theorem}[\cite{blackwell1950comparison}]\label{thm:blackwell}
Let $P$ and $Q$ be probability distributions on some space and $P'$ and $Q'$ be probability distributions on another space. The two statements below are equivalent to each other:
\begin{enumerate}
\item[(a)] $T(P,Q) \le_{\succeq} T(P',Q')$.
\item[(b)] There exists a (randomized) algorithm $\mathrm{Proc}$ such that $\mathrm{Proc}(P)=P', \mathrm{Proc}(Q)=Q'$.
\end{enumerate}
\end{theorem}

This theorem is a special case of a more general theorem called Blackwell's informativeness theorem, which concerns the ranking of information structures and establishes an equivalence between rankings from different perspectives \citep{blackwell1950comparison}. The particular version presented in Theorem \ref{thm:blackwell} is proven in \cite{blackwell1979theory} (Chapter 12, page 334). Its connection to differential privacy was first recognized by \cite{KOV} and further developed by \cite{dong2022gaussian}.

Now we briefly explain how this theorem implies the representation theorem. See Proposition~B.1 of \citet{dong2022gaussian} for more details.

Consider two hypothesis testing problems,\footnote{Either with two algorithms or two different pairs of neighboring datasets.} with probability distributions $(P_1, Q_1)$ and $(P_2, Q_2)$, respectively. A contradiction to Theorem~\ref{thm:full} would occur only if $T(P_1, Q_1) = T(P_2, Q_2)$ while $D(P_1, Q_1) \ne D(P_2, Q_2)$.

However, this is impossible. If $T(P_1, Q_1) = T(P_2, Q_2)$, which gives $T(P_1, Q_1) \le_\succeq T(P_2, Q_2)$, Blackwell's theorem shows that there must exist $\mathrm{Proc}$ such that $\mathrm{Proc}(P_1)=P_2, \mathrm{Proc}(Q_1)=Q_2$. By Axiom~\ref{axim}, the metric $D$ used in the differential privacy definition under the consideration of Theorem~\ref{thm:full} must satisfy
\[
D(P_1, Q_1) \succeq D(\mathrm{Proc}(P_1), \mathrm{Proc}(Q_1)) = D(P_2, Q_2)
\]
due to the post-processing property. Likewise, we have
\[
D(P_1, Q_1) \preceq  D(P_2, Q_2).
\]
Taken together, the two latest displays show that $D(P_1, Q_1) =  D(P_2, Q_2)$ if $T(P_1, Q_1) = T(P_2, Q_2)$. This proves the representation theorem.

Blackwell's theorem shows that the trade-off function induces the same ordering as post-processing. As a direct consequence, $T(P,Q) \le_\succeq T(P',Q')$ implies $D(P,Q) \succeq D(P',Q')$ under Axioms~\ref{axiom:ht} and \ref{axim}. In sharp contrast, this is not true for R\'enyi divergences, even taken over all orders $\gamma > 1$ collectively.
\begin{proposition}[\cite{concentrated2,dong2022gaussian}]\label{prop:renyi_fail}
Let $P_\ep$ and $Q_{\ep}$ be Bernoulli distributions that take value 1 with probabilities $\frac{\e^\ep}{1+\e^\ep}$ and $\frac{1}{1+\e^\ep}$, respectively. There exists $\ep>0$ such that the two statements below both hold:
\begin{enumerate}
\item[(a)]
For all $\gamma>1$, $D_\gamma\big(\N(0,1), \N(\ep,1)\big) \ge D_\gamma(P_\ep, Q_\ep)$.
\item[(b)] For total variation distance, $\TV\big(\N(0,1), \N(\ep,1)\big) < \TV(P_\ep, Q_\ep)$.
\end{enumerate}
\end{proposition}


\section{Toolbox for $f$-Differential Privacy}
\label{sec:privacy-accounting}

In this section, we review basic statistical properties of $f$-differential privacy and an important desideratum of this privacy framework under composition.

\subsection{Basic Properties}
\label{sec:basics-f-dp}

As the name suggests, an alternative definition of a trade-off function shows that it gives the optimal trade-off between the type I and type II errors for testing the hypothesis testing problem in Equation~\ref{eq:ht}. Consider a rejection rule $0 \le \phi \le 1$ that takes as input the output of the algorithm $M$, with its type I and type II errors defined as
\[
\alpha_\phi = \E_{P} \phi, \quad \beta_\phi = 1 - \E_{Q} \phi,
\]
respectively. The trade-off function of $P$ and $Q$ can be defined as
\[
T(P, Q)(\alpha) = \inf_{\phi} \left\{ \beta_\phi: \alpha_\phi \leqslant \alpha \right\}
\]
for any $0 \le \alpha \le 1$. By the Neyman--Pearson lemma, this definition is equivalent to the likelihood-based definition in Equation~\ref{eq:lration}. If $T(P,Q) \le_\succeq T(P', Q')$, then in a very strong sense, $P$ and $Q$ are easier to distinguish than $P'$ and $Q'$ at \textit{any} level of type I error. In general, two trade-off functions may not be comparable since $\le_\succeq$ or, equivalently, $\ge_\preceq$ is a partial order.

Let $f$ be a trade-off function, which is convex, continuous, non-increasing, and satisfies $f(\alpha) \le 1 - \alpha$ for all $\alpha$. Below, we give the definition of $f$-differential privacy.
\begin{definition}[\cite{dong2022gaussian}] \label{def:kstable}
An algorithm $M$ is said to be $f$-differentially private if
\begin{equation}\label{eq:tf1}
T\big(M(S), M(S')\big) \ge_\preceq f
\end{equation}
for any neighboring datasets $S$ and $S'$.
\end{definition}

\begin{figure}[!htp]
\centering
\includegraphics[scale=0.6]{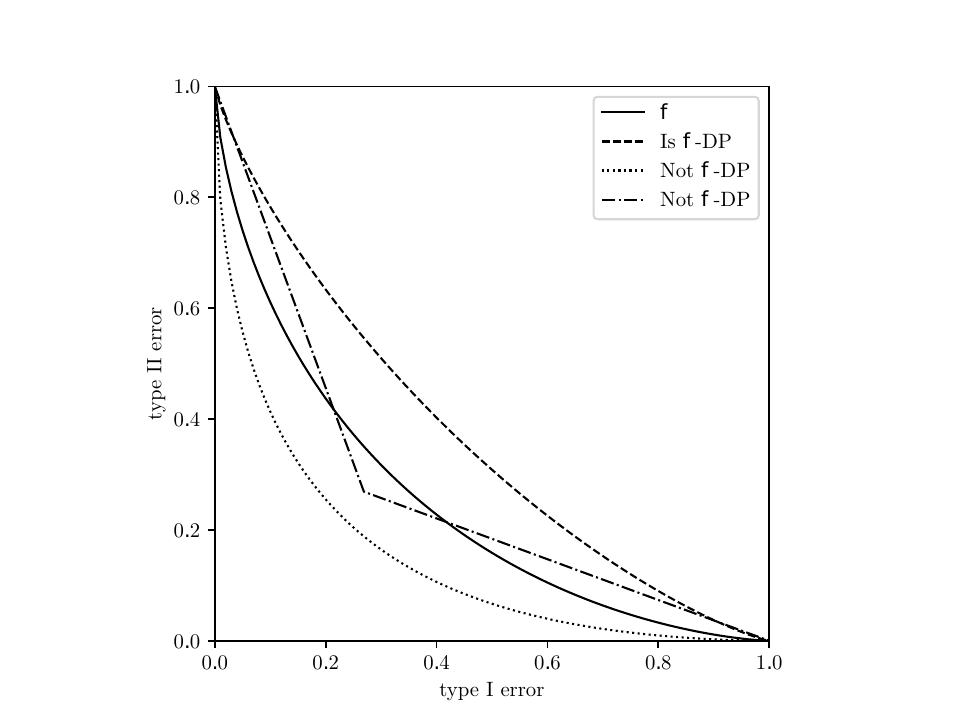}
\caption{The graph shows three different trade-off functions for $T(M(S), M(S'))$. Among these, only the dashed line represents a trade-off function that satisfies $f$-differential privacy. Adapted from Figure 2 in \cite{dong2022gaussian}.}
\label{fig:nonfdp}
\end{figure}

Figure~\ref{fig:nonfdp} gives an illustration of this definition.

Writing the inverse $f^{-1} = T(Q, P)$ if $f = T(P, Q)$, or equivalently, $f^{-1}(\beta) = \alpha$ if $f(\alpha) = \beta$. In general, we can assume that $f$ is symmetric in the sense that if $f = f^{-1}$ due to the symmetry of $S$ and $S'$. To see this point, note that Definition~\ref{def:kstable} shows that $T\big(M(S'), M(S)\big) \ge_\preceq f$, which gives
\begin{equation}\label{eq:tf2}
T\big(M(S), M(S')\big) \ge_\preceq f^{-1}
\end{equation}
by recognizing the simple fact that $f_1 \ge_\preceq f_2$ implies $f_1^{-1} \ge_\preceq f_2^{-1}$. Taken together, Equations~\ref{eq:tf1} and \ref{eq:tf2} show that 
\[
T\big(M(S), M(S')\big) \ge_\preceq \max\{f, f^{-1}\}.
\]
Above $\max\{f, f^{-1}\}$ is also a trade-off function~\citep{dong2022gaussian} and yields a tighter bound since $\max\{f, f^{-1}\} \ge_\preceq f$.

As discussed in Section~\ref{sec:metrics}, $f$-differential privacy implies any other differential definitions, including $(\epsilon, \delta)$-differential privacy via Equation~\ref{eq:dual}. We ask what the least strong $f$-differential privacy guarantee is that can lead to $(\epsilon, \delta)$-differential privacy for some given $\epsilon$ and $\delta$. We have the following result.
\begin{theorem}[\cite{wasserman_zhou,KOV,dong2022gaussian}]
Fix $\ep \ge 0$ and $0 < \delta \le 1$. If $f$-differential privacy implies $(\epsilon, \delta)$-differential privacy in the sense that $\epsilon$ is the root to Equation~\ref{eq:dual}, then it must hold that
\[
f(\alpha) \ge f_{\ep, \delta}(\alpha) := \max\left\{ 0,1 - \delta - \e^\ep \alpha, \e^{-\ep}(1-\delta-\alpha) \right\}
\]
for any $0 \le \alpha \le 1$.
\end{theorem}

This bound is tight because Equation~\ref{eq:dual} is satisfied by $f_{\ep, \delta}$. Indeed, let $P_{\ep, \delta}$ and $Q_{\ep, \delta}$ be such that $f_{\ep, \delta} = T(P_{\ep, \delta}, Q_{\ep, \delta})$. Then under either $P_{\ep, \delta}$ or $Q_{\ep, \delta}$, the likelihood ratio between $P_{\ep, \delta}$ and $Q_{\ep, \delta}$, other than an event of probability $\delta$ where they are singular, must be either $\e^{\ep}$ or $\e^{-\ep}$ almost surely \citep{wasserman_zhou}. This shows that $(\ep,\delta)$-differential privacy is equivalent to $f_{\ep, \delta}$-differential privacy in a rigorous sense, demonstrating that $f$-differential privacy is an extension of $(\ep, \del)$-differential privacy.

\subsection{Achieving Differential Privacy}\label{sec:achi-diff-priv}

A template for achieving differential privacy is to add additive noise to the output of the algorithm, which leads to the class of additive mechanisms.\footnote{It is worthwhile mentioning that there are important non-additive mechanisms such as the exponential mechanism~\citep{mcsherry2007mechanism}.} For instance, $\epsilon$-differential privacy can be achieved by using Laplace noise. Another example is the Gaussian mechanism, which adds Gaussian noise, distorting the output of the algorithm significantly less in the rare event because Gaussian noise has a lighter tail than Laplace noise. Notably, the Gaussian mechanism is not $\epsilon$-differential privacy for any $\ep < \infty$, which motivates $(\ep, \del)$-differential privacy in the history of the development of differential privacy~\citep{approxdp}.

Let the algorithm $M(S)$ take the form $M(S) = \widebar{M}(S) + Z$, where $\widebar{M}$ is a statistic of interest, say, the average height of individuals in $S$, and $Z$ is some additive noise. The privacy guarantee of $M$ depends on the sensitivity of $\widebar{M}$:
\[
\mathrm{sens}(\widebar{M}) = \sup_{S, S'} |\widebar{M}(S) - \widebar{M}(S') |,
\]
where the supremum is taken over neighboring $S$ and $S'$ and $|\cdot|$ is a norm. Taking $Z \sim \N(0, \mathrm{sens}(\widebar{M})^2/\mu^2)$ for some $\mu > 0$, the Gaussian mechanism $M$ satisfies
\[
T(M(S), M(S')) \ge_\preceq T(\N(0, 1), \N(\mu, 1)) := G_{\mu}
\]
for any neighboring $S$ and $S'$. Above, $G_\mu(\alpha) = \Phi\big(\Phi^{-1}(1-\alpha)-\mu\big)$ is called the Gaussian trade-off function, where $\Phi$ denotes the standard normal cumulative distribution function (CDF). Thus, this Gaussian mechanism is $G_{\mu}$-differentially private.

We call an algorithm $\mu$-Gaussian differentially private if it is $G_{\mu}$-differentially private. This yields an important family of $f$-differential privacy guarantees. It offers the tightest possible characterization of the Gaussian mechanism's privacy guarantees. Its privacy properties are fully captured by a single parameter---the mean of a unit-variance Gaussian distribution. This concise representation facilitates easy interpretation and communication of the privacy level ensured by Gaussian differential privacy.

Given a trade-off function $f$, an important question is finding a noise distribution $Z$ such that the output 
\[
M(S) = \widebar{M}(S) + \mathrm{sens}(\widebar{M}) Z
\]
is $f$-differentially private. \citet{awan2023canonical} proved the existence of such distributions for a general trade-off function $f$ that is symmetric and nontrivial---that is, there exists some $\alpha$ such that $f(\alpha) < 1 - \alpha$. Specifically, such a distribution satisfies:
\begin{itemize}
\item[1.] $T(Z, Z + \Delta) \ge_\preceq f$ for any $0 \le \Delta \le 1$, with equality at $\Delta = 1$, 
\item[2.] $T(Z, Z+1)(\alpha) = F(F^{-1}(\alpha) - 1)$, where $F$ is the CDF of $Z$, and 
\item[3.] $Z$ is symmetric with respect to 0. 
\end{itemize}

In \citet{awan2023canonical}, the authors called $Z$ a canonical noise distribution for $f$ if the three properties above are satisfied. Notably, a canonical noise distribution fully exploits the privacy budget conveyed by the trade-off function. This distribution admits an explicit construction based on the given trade-off function \citep{awan2023canonical}. It has been extended to multivariate settings and shown to attain certain optimal statistical properties \citep{awan2022log,awan2023optimizing}.

\subsection{Composition}
The composition of multiple privacy-preserving analyses on the same dataset leads to a degradation of the overall privacy guarantee, as intermediate results are released. In the $f$-differential privacy framework, the composition of algorithms can be formalized using the tensor product operation $\otimes$ on the associated trade-off functions. Specifically, if algorithm $M_1$ satisfies $f$-differential privacy with $f = T(P, Q)$ and algorithm $M_2$ satisfies $g$-differential privacy with $g = T(P', Q')$, then their composition satisfies $(f \otimes g)$-differential privacy, where~\citep{dong2022gaussian}
\[
f \otimes g := T(P \times P', Q \times Q').
\]

More generally, let $M_i$ be an $f_i$-differentially private algorithm conditionally on any output of the prior algorithms for $i = 1, \ldots, m$. Then the $m$-fold composition of these algorithms satisfies $f_1 \otimes \cdots \otimes f_m$-differential privacy. The tensor product operation $\otimes$ is commutative and associative. This privacy bound is tight in general, in contrast to $(\ep, \del)$-differential privacy, where the composition mechanism cannot be tightly characterized by a single pair of privacy parameters \citep{KOV}. For the special case of Gaussian differential privacy, composition takes a simple form: 
\[
G_{\mu_1} \otimes \cdots \otimes G_{\mu_m} = G_{\sqrt{\mu_1^2 + \cdots + \mu_m^2}}.
\]

Importantly, a central limit theorem phenomenon emerges when composing many ``very private'' $f$-differentially private algorithms. Let 
\[
\kl(f) = -\int_0^1 \log(-f'(x)) \d x \text{ and } \kappa_2(f) = \int_0^1 \log^2(-f'(x)) \d x.
\] 
Informally, if each trade-off function $f_i$ is close to perfect privacy---that is, $\id(\alpha) = 1-\alpha$ for all $0 \le \alpha \le 1$---then for a sufficiently large number of compositions $m$, we have the approximation \citep{dong2022gaussian}:
\[
f_1 \otimes f_2 \otimes \cdots \otimes f_m \approx G_\mu,
\]
where
\[
\mu \approx \frac{2 \sum_{i=1}^m \kl(f_i)}{\sqrt{\sum_{i=1}^m \kappa_2(f_i)}}.
\]
This approximation is particularly relevant for privacy analysis in deep learning, where training typically involves tens of thousands of iterations.

To improve the accounting of privacy losses under composition, both analytical approaches based on characteristic functions and numerical approaches using the fast Fourier transform have been developed~\citep{koskela2020computing,gopi2021numerical,zhu2022optimal}. The Edgeworth accountant provides a way to combine these analytical and numerical benefits by first expressing the overall privacy loss in tensor product form and then tightly converting it to $(\ep,\delta)$-differential privacy guarantees.

\begin{proposition}[\cite{wang2022analytical}]\label{prop:eps_delta_f_dual}
Let $F_{X, m}, F_{Y, m}$ be the CDFs of $X_{1}+\cdots+X_{m}$ and $Y_{1}+\cdots+Y_{m}$, respectively.
Then, $ f_1 \otimes \cdots \otimes f_m$-differential privacy implies $(\epsilon, \delta)$-differential privacy for all $\epsilon> 0$ with
\begin{equation}\nonumber
\delta = 1 - F_{Y,m}(\epsilon) - \e^{\epsilon}(1-F_{X,m}(\epsilon)).
\end{equation}
\end{proposition}

\begin{figure}[!htp]
\centering
\includegraphics[scale=0.65]{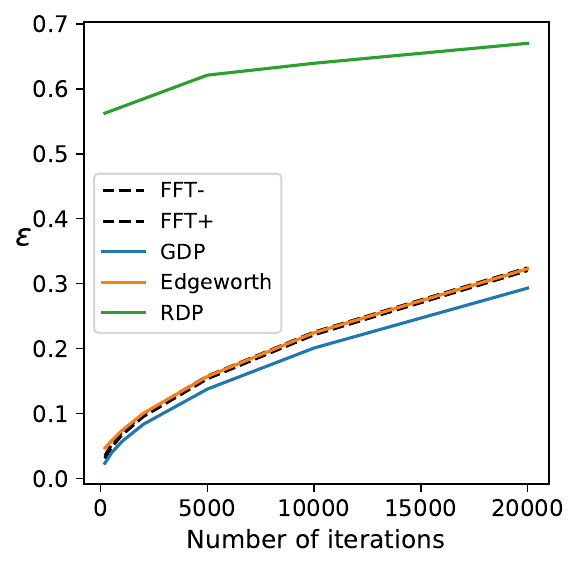}
\caption{Comparison between different privacy accountants in terms of the privacy parameter $\epsilon$ with $\delta = 10^{-5}$ in private federated analytics. GDP, Edgeworth, and RDP correspond to the privacy bounds obtained from \citet{bu2020deep}, \citet{wang2022analytical}, and \citet{wang2018subsampled}, respectively. FFT- and FFT+ denote the lower and upper bounds derived by numerical methods \citep{koskela2020computing,gopi2021numerical}, which sandwich the Edgeworth bound tightly. For experimental details, see \citet{wang2022analytical}.}
\label{fig:edgeworth}
\end{figure}

This result reduces the computation of $(\ep, \del)$-differential privacy guarantees to the evaluation of the CDFs, $F_{X,m}$ and $F_{Y,m}$. \citet{wang2022analytical} proposed approximating these CDFs using first-order Edgeworth expansions. An appealing feature of this approach is that a finite-sample error bound is derived, which is of order $O(1/m)$ when the trade-off functions $f_i$ are the same, often the case in applications. Although the expression formula is somewhat complex, another advantage is that the Edgeworth accountant is analytic, so the computational cost is nearly zero. The performance of this privacy accountant is illustrated in Figure~\ref{fig:edgeworth}, where all privacy bounds are converted to $(\ep, \del)$-differential privacy.


\section{Applications of $f$-Differential Privacy}
\label{sec:applications}

In this section, we showcase several applications where the use of $f$-differential privacy enhances privacy bounds compared to earlier approaches. The focus is to glean the statistical structures of the problems that fit into the scope of $f$-differential privacy. The elaboration of each application is brief, and we refer interested readers to the references for more details.

\subsection{Private Deep Learning}\label{sec:deep-learning}

Due to their computational scalability, stochastic gradient descent (SGD) and its variants, such as Adam \citep{kingma2014adam}, are used to minimize the empirical risk
\[
L(\theta) = \frac1n \sum_{i=1}^n \ell(\theta, x_i)
\]
over the weights $\theta$ in some Euclidean space for training deep neural networks, where $\ell$ is a loss function and $\{x_1, \ldots, x_n\}$ denotes the dataset $S$. In each iteration $t$ of SGD, a mini-batch $I_t$ is randomly sampled from $\{1, 2, \ldots, n\}$ with subsampling probability $p$, and SGD updates the weights according to
\[
\theta_{t+1} = \theta_t - \frac{\eta_t}{|I_t|} \sum_{i \in I_t} \nabla_{\theta} \ell(\theta_t, x_i),
\]
where $\eta_t$ denotes the learning rate, $|I_t|$ is the mini-batch size, and $\nabla_{\theta} \ell$ denotes the subgradient of $\ell$ with respect to $\theta$.

To achieve differential privacy, the per-example gradients are clipped to have $\ell_2$ norm no more than $R$, which ensures a finite sensitivity, and Gaussian noise of scale $\sigma R$ is added to the average clipped gradient, which is equivalent to the Gaussian mechanism \citep{chaudhuri2011differentially,song2013stochastic,bassily2014private,deep}. The iterative nature of this variant of SGD, referred to as DP-SGD, renders the privacy analysis in the natural form of composition. Under the $f$-differential privacy framework, writing $p = B/n$ in the case of a fixed mini-batch size $|I_t| \equiv B$, \citet{dong2022gaussian} and \citet{bu2020deep} showed that running DP-SGD for $T$ iterations is $\min\{f,f^{-1}\}^{\ast\ast}$-differentially private, where 
\[
f = \big(pG_{1/\sigma}+(1-p)\id\big)^{\otimes T}.
\]
Above, recall that $\id(\alpha) = 1-\alpha$ for $0 \le \alpha \le 1$.

To evaluate this privacy bound, we can derive an approximation using the privacy central limit theorem in an asymptotic regime where $p \sqrt{T} \to \nu$ for some constant $\nu > 0$ as $T \to \infty$. In this setting, $\big(pG_{1/\sigma}+(1-p)\id\big)^{\otimes T}$ converges to $G_{\mu}$ for 
\[
\mu = \nu\sqrt{\e^{1/\sigma^2}-1}.
\] 
Consequently, the overall privacy guarantee simplifies to 
\[
\min\{f,f^{-1}\}^{\ast\ast} \approx \min\{G_{\mu}, G_{\mu}^{-1}\}^{\ast\ast} = G_{\mu}.
\] 
This shows that private deep learning using DP-SGD is approximately $\frac{B}{n}\sqrt{T(\e^{1/\sigma^2}-1)}$-Gaussian differentially private.

In \citet{wang2022analytical}, an improved approximation of $\min\{f,f^{-1}\}^{\ast\ast}$ is obtained using the Edgeworth accountant. This is done by approximating $f = \big(pG_{1/\sigma}+(1-p) \id\big)^{\otimes m}$ and $f^{-1} = \{\big(pG_{1/\sigma}+(1-p) \id\big)^{-1}\}^{\otimes T}$ separately. For $f$, the privacy-loss log-likelihood ratios are 
\[
\begin{aligned}
&X_i^{(1)} = \log(1-p+p\e^{\frac{\xi_i}{\sigma}-\frac{1}{2\sigma^2}})\\
&Y_i^{(1)} = \log(1-p+p\e^{\frac{\zeta_i}{\sigma}-\frac{1}{2\sigma^2}})\\
\end{aligned}
\] 
for $1\le i \le m$, with independent and identically distributed (i.i.d.) $\xi_i\sim \mathcal{N}(0, 1)$ and $\zeta_i\sim (1-p)\mathcal{N}(0, 1) + p \mathcal{N}(1/\sigma,1)$. For $f^{-1}$, the privacy-loss log-likelihood ratios are $X_i^{(2)} = - Y_i^{(1)}$ and $Y_i^{(2)} = -X_i^{(1)}$. The remaining step is to evaluate and approximate the moment-generating functions of these privacy-loss log-likelihood ratios. 

An illustration in Figure~\ref{fig:edgeworth} shows how the Edgeworth accountant improves on the central limit theorem-based approximation, and both of these are far superior to the accounting approach based on R\'enyi differential privacy~\citep{wang2019subsampled}.

\subsection{Private Convex Optimization}
\label{sec:last-iterate-with}

In Section \ref{sec:deep-learning}, the privacy bound in Gaussian differential privacy satisfies
\[
\frac{B}{n}\sqrt{T(\e^{1/\sigma^2}-1)} = O(\sqrt{T}) \goto \infty
\] 
as the number of iterations $T \goto \infty$ while any other parameters are fixed, meaning that the privacy guarantee of DP-SGD becomes vacuous when $T$ is very large. This is essentially due to the release of all intermediate iterates so that the privacy loss always accumulates during iteration. For strongly convex optimization, however, it has been shown that releasing the last iterate has bounded privacy loss due to a certain contraction effect, thereby leading to significant privacy amplification \citep{feldman2018privacy,altschuler2022privacy}. Their techniques are based on R\'enyi differential privacy.

In a recent work, \citet{bok2024shifted} used $f$-differential privacy to improve the privacy analysis of the last iterate in minimizing a strongly convex function. The core technique of \citet{bok2024shifted} is to construct an auxiliary process $\{ \widetilde{\theta}_t \}$ that interpolates between the two sequences of iterates $\{ \theta_t \}$ and $\{ \theta_t' \}$, which correspond to the null and the alternative hypotheses. This auxiliary process is designed to match the second sequence until an intermediate time point $\tau$ and coincides with the first sequence at the final time $T$. The key advantage of this approach is that it allows us to unroll the argument backward from $T$ to $\tau$, instead of all the way to the initial point.

This technique effectively replaces the divergent $\sqrt{T}$ term found in previous $f$-differential privacy bounds with a term that scales with $T - \tau$. The choice of the intermediate time point $\tau$ is a tunable parameter in the analysis. On the one hand, selecting a larger value of $\tau$ reduces the number of steps one needs to unroll. On the other hand, a smaller $\tau$ provides the auxiliary process $\{\widetilde{\theta}_t \}$ with more iterations to gradually transition from $\theta_{\tau}'$ to $\theta_T$ with more flexibility.

\begin{table}[H]
\centering
\begin{tabular}{llll}
\toprule
Epochs & \multicolumn{1}{c}{Composition (all iterates)} & \multicolumn{1}{c}{R\'enyi differential privacy} & \multicolumn{1}{c}{$f$-differential privacy} \\
 \midrule
50 & 30.51 & 5.82 & 4.34 \\
100 & 49.88 & 7.61 & 5.60 \\
200 & 83.83 & 9.88 & 7.58 \\
\bottomrule
\end{tabular}
\vskip.5em
\caption{Privacy bounds, measured by $\ep$ in $(\epsilon, \delta)$-differential privacy with $\delta = 10^{-5}$, for a cyclic version of noisy gradient descent applied to regularized logistic regression. Note that one epoch is a complete pass through the entire training dataset. The $f$-differential privacy based approach yields stronger privacy guarantees compared to the other two bounds, which are based on composition and R\'enyi differential privacy, respectively. See more details on the experimental setup in \cite{bok2024shifted}.}
\label{tab:1}
\end{table}

Table~\ref{tab:1} shows the results from an experiment of regularized logistic regression on MNIST \citep{lecun1998mnist}. By converting privacy guarantees to $(\ep, \del)$-differential privacy using the duality, the results show that this $f$-differential privacy-powered technique obtains tighter privacy bounds than those obtained using R\'enyi differential privacy and, consequently, enables longer training time with the same privacy budget.

\subsection{Shuffled Mechanism}
\label{sec:mixture-distribution}

In many applications, the output of the algorithm has a mixture probability distribution. Formally, write $P = \sum_{i=1}^k w_i P_i$ and $Q = \sum_{i=1}^k w_i Q_i$ for some weights $w_i \geq 0$ such that $\sum_{i=1}^k w_i = 1$ and some simple distributions $P_i$ and $Q_i$. \citet{wang2023unified} obtained an expression that relates the trade-off function of the mixture distributions to those of the mixture components:
\begin{align}\label{eq:mixture}
T(P,Q)(\alpha(r,c)) \ge_\preceq \sum_{i=1}^k w_i T(P_i,Q_i)(\alpha_i(r,c))
\end{align}
for $r \ge 0$ and $0 \le c \le 1$, where 
\[
\alpha_i(r,c) = \P_{P_i}\left( \frac{\d Q_i}{\d P_i} > r\right) + c \P_{P_i}\left(\frac{\d Q_i}{\d P_i} = r \right)
\]
is the type I error using the likelihood ratio for testing $P_i$ versus $Q_i$ and $\alpha (r,c) = \sum_{i=1}^k w_i \alpha_i(r,c)$. The parameter $c$ in the definition is redundant when the Radon--Nikodym derivative $\frac{\d Q_i}{\d P_i}$ does not have any point mass. This bound is tight in the sense that the inequality in Equation~\ref{eq:mixture} can reduce to an equality for certain instances of mixture distributions \citep{wang2023unified}.

\citet{wang2023unified} applied this tight bound under the $f$-differential privacy framework to analyze shuffled mechanisms \citep{erlingsson2019amplification}. Roughly speaking, a shuffled mechanism is a simple algorithmic idea that randomly permutes a collection of locally differentially private reports to boost privacy guarantees, and the shuffling process naturally endows the output with a structure of mixture distributions~\citep{feldman2022hiding,feldman2023stronger}.

A comparison with \citet{feldman2023stronger} is performed, with results presented in Tables \ref{tab:compareFeldman_10k} and \ref{tab:compareFeldman_numerical}. As earlier in Section~\ref{sec:last-iterate-with}, the privacy bounds have been converted to $(\ep, \del)$-differential privacy using the dual representation in Equation~\ref{eq:dual}. The results show that $f$-differential privacy analysis yields significantly tighter privacy bounds than R\'enyi differential privacy.

\begin{table}[H]
  \centering
  \begin{tabular}{ lllllll  }
\toprule 
 $\epsilon$  & 0.5 & 0.6 & 0.7 &0.8 & 0.9 & 1.0 \\
 \hline
 $\delta$ (R\'enyi differential privacy)  & $0.9494$ & $0.3764$ & $0.1038$  & $0.0181$ & $0.0018$  & $8 \times 10^{-5}$ \\
 \hline
 $\delta$ ($f$-differential privacy) & $3 \times 10^{-6}$ & $10^{-7}$ & $4 \times 10^{-9}$ & $9 \times 10^{-11}$ & $2 \times 10^{-12}$ & $2 \times 10^{-14}$  \\
  \bottomrule
\end{tabular}
\vskip.5em
\caption{Comparison between the privacy bounds of \cite{feldman2023stronger} and \cite{wang2023unified} on shuffled mechanisms, in terms of $\delta$ for the same $\ep$. See details on the experimental setup in \cite{feldman2023stronger}.}
\label{tab:compareFeldman_10k}
\end{table}

\begin{table}[H]
  \centering
  \begin{tabular}{ llllll  }
\toprule 
 $\delta$ & $5 \times 10^{-5}$& $3 \times 10^{-6}$ & $10^{-7}$ & $4 \times 10^{-9}$ & $9 \times 10^{-11}$   \\
 \hline
  $\epsilon$ ($f$-differential privacy) & $0.4$ & 0.5 & 0.6 & 0.7 &0.8 \\
 \hline
 Numerical $\epsilon$ upper bound & $1.014$& $1.085$ & $4.444$ & $4.444$  & $4.444$ \\
  \hline
 Numerical $\epsilon$ lower bound & $0.369$& $0.470$ & $0.575$ & $0.664$  & $0.758$\\
  \bottomrule
\end{tabular}
\vskip.5em
\caption{Comparison between the privacy bounds of \cite{feldman2023stronger} and \cite{wang2023unified} on shuffled mechanisms, in terms of $\ep$ for the same $\del$. The numerical upper bound on $\ep$ is from \cite{feldman2023stronger} using analysis of R\'enyi differential privacy, whereas the lower bound is obtained by binary search. When $\delta \le 10^{-7}$, the upper bound of \cite{feldman2023stronger} becomes invalid, and is therefore replaced by a trivial bound of value $4.444$. In contrast, the bound of \cite{wang2023unified} is valid for any value of $\delta$.}
\label{tab:compareFeldman_numerical}
\end{table}

\subsection{U.S.\ Census Data}
\label{sec:us-census}

The Census Bureau collects vast amounts of personal data and has a responsibility to protect the privacy and confidentiality of the individual responses contained within its records. The Census is susceptible to potential attacks that could compromise this sensitive information \citep{abowd2018us}. To address this vulnerability and bolster privacy safeguards, the Census Bureau has adopted differential privacy within its disclosure avoidance system for enhancing the privacy protections surrounding the 2020 Census data.

In \citet{su2024revealing}, an improved privacy analysis in $f$-differential privacy is obtained. The improved privacy analysis is due to the released census naturally having a composition structure. In the disclosure avoidance system (DAS) employed by the Census Bureau \citep{abowd20222020}, each query $M(S)$ is privatized by adding discrete Gaussian noise, resulting in $M(S) = \widebar{M}(S) + \mathcal{N}_{\mathbb{Z}}(0,\sigma^2)$. Here, a discrete Gaussian distribution $\mathcal{N}_{\mathbb{Z}}(\mu, \sigma^2)$ with mean $\mu$ and proxy variance $\sigma^2$ has probability mass function \citep{canonne2020discrete}
 \begin{align*}
    p_{\mu,\sigma}(x) = \frac{\e^{-(x-\mu)^2/2\sigma^2}}{{\sum_{i\in\mathbb{Z}}\e^{-(i-\mu)^2/2\sigma^2}}}
\end{align*}
for any integer $x$.

The DAS utilizes a sequence of counting queries $\left\{\widebar{M}_i(S) \right\}_{i=1}^m$ with $m=72$, as specified in the privacy allocation file \citep{privacyallocation2022}. These queries cover 8 geographical levels, including nation, state, and county, with each level containing 9 queries. Each query is privatized as $M_i(S) = \widebar{M}_i(S) + \mathcal{N}_{\mathbb{Z}}(0,\sigma_i^2)$. As each $\widebar{M}(S)$ is a counting query with sensitivity 1, the analysis of the privacy bound boils down to evaluating the trade-off function for testing the following hypotheses:
\begin{align*}
H_0: P = \bigotimes_{i=1}^m \mathcal{N}_{\mathbb{Z}}(0,\sigma_i^2) \qquad \text{vs.} \qquad H_1: Q = \bigotimes_{i=1}^m \mathcal{N}_{\mathbb{Z}}(1,\sigma_i^2).
\end{align*}

Several challenges arise when evaluating the trade-off function. First, the likelihood ratio involves a quantity that no longer follows a discrete Gaussian distribution, in contrast to the Gaussian case where any linear combination of Gaussian variables remains Gaussian. Moreover, except in the i.i.d.\ case where $\sigma_i \equiv \sigma$, the support of that quantity is a complex grid. \citet{su2024revealing} introduced several new techniques to circumvent these challenges.

In their experiments, \citet{su2024revealing} showed that this $f$-differential privacy-based approach significantly outperforms the R\'enyi differential privacy-based approach \citep{canonne2020discrete} when both privacy guarantees are converted to $(\ep, \del)$-differential privacy. This, in turn, enables one to add less noise to the raw census data. The results on the noisy measurement files demonstrate about a 10\% reduction in the mean squared error when $f$-differential privacy is used in place of R\'enyi differential privacy~\citep{su2024revealing}. As an example, the privacy guarantee for the US level geography on the noisy measurement files is improved from $(2.63, 10^{-11})$-differential privacy to $(2.32, 10^{-11})$-differential privacy.


\section{Concluding Remarks}
\label{sec:future-directions}

In this review paper, we have explored differential privacy from a statistical viewpoint, motivated by the hypothesis testing framework. Drawing primarily from the work of \citet{wasserman_zhou,KOV,dong2022gaussian}, we have shown that hypothesis testing is not merely a convenient tool but in fact the right language for reasoning about differential privacy. In light of this, $f$-differential privacy becomes a very natural definition that preserves both the hypothesis testing interpretation and the property of post-processing. Owing to its coherent mathematical structure, the $f$-differential privacy framework is accompanied by powerful tools, such as composition theorems, which enable tight and even lossless privacy analysis. This paper has demonstrated the effectiveness of this privacy definition in several important applications. Nevertheless, the review does not cover other notable applications, including federated learning, privacy auditing, and Markov chain Monte Carlo \citep{zheng2021federated,nasr2023tight,lin2023tractable}.

In closing, we discuss several future research directions. First, while semantically being the most informative, $f$-differential privacy comes with increased technical difficulty in evaluating trade-off functions compared to divergence-based differential privacy definitions. Indeed, \citet{bok2024shifted} and \citet{su2024revealing} (discussed in Section~\ref{sec:applications}) devoted much effort to tackling this difficulty. Therefore, it would be interesting to develop tools that facilitate the use of $f$-differential privacy. For example, more efficient numerical methods for composing privacy bounds would be welcome. Second, $f$-differential privacy often leads to asymptotic privacy bounds \citep{bu2020deep,wang2022analytical}. An interesting fact is that asymptotic bounds are often sharper than non-asymptotic bounds, which leads to a trade-off between the usefulness and interpretability of the privacy bounds. An important research question is to study to what extent this trade-off is fundamental. Last, it would be interesting to extend the $f$-differential privacy framework to non-Euclidean data \citep{jiang2023gaussian,liu2022identification}. All in all, by adapting and refining the concepts and techniques developed for $f$-differential privacy, researchers can potentially unlock new possibilities for privacy-preserving data analysis across a wide range of domains and applications.

\section*{Acknowledgments}

This work was supported in part by NSF DMS-2310679, a Meta Faculty Research Award, and Wharton AI for Business. The author would like to thank Gautam Kamath, Adam Smith, and Chendi Wang for helpful comments on an early version of the paper.

{\small
\bibliographystyle{alpha}
\bibliography{ref}
}

\end{document}